\documentclass[a4paper,12pt]{article}
\usepackage{amsmath,amsfonts,amssymb}
\usepackage[english,russian]{babel}
\usepackage{graphicx}

\makeindex

\begin{document}
\sloppy

\begin{center}

\textbf{The modified problems for the equation of Euler--Darboux\\
in the case of parameters on the module equal to 1/2}\\
\medskip

\textbf{Dolgopolov M.V.$^{1,2}$, Rodionova I.N.$^2$}
\medskip

$^1$\textit{Samara State Technical University}, \texttt{mikhaildolgopolov68@gmail.com}\\
$^2$\textit{Research Laboratory of Mathematical Physics}, \texttt{mvdolg@yandex.ru}\\
\end{center}

\vskip10pt

\textbf{Abstract}. We consider the Euler--Darboux equation with parameters modulo $\displaystyle\frac{1}{2}$ and generalization to the space 3D analogue. Due to the fact that the Cauchy problem in its classical formulation is incorrect for such parameter values, the authors propose formulations and solutions of modified Cauchy-type problems with parameter values: a) $\alpha=\beta=\displaystyle\frac{1}{2}$, b) $\alpha=-\,\displaystyle\frac{1}{2}$, $\beta=+\,\displaystyle\frac{1}{2}$, c)~$\alpha=\beta=-\,\displaystyle\frac{1}{2}$.
The obtained result is used to formulate an analogue of the $\Delta_1$ problem in the first quadrant with the setting of boundary conditions with displacement on the coordinate axes and non-standard conjugation conditions on the singularity line of the coefficients of the equation $y=x$. The~first of these conditions glues the normal derivatives of the desired solution, the~second contains the limit values of the combination of the solution and its normal derivatives. The problem was reduced to a uniquely solvable system of integral equations.\\


%
%
%
%
%
%
%
%
%
%
%

\textbf{Introduction}. Euler--Poisson--Darboux equations
\begin{equation}
U_{x y} + \,\frac{\beta}{y- x} U_x - \, \frac{\alpha}{y- x} U_y  =0 \label{eq:1}
\end{equation}
are widely used in gas and hydrodynamics, shell theory, in various sections of continuum mechanics.

Due to the fact that degenerate equations of hyperbolic type in characteristic coordinates are reduced to the equation (\ref{eq:1}), many Soviet and foreign mathematicians were engaged in the study of boundary value problems for the Euler-Darboux equation. A~detailed bibliography on this subject is contained in the monograph of M.\,M.\,Smirnov~[1]. The paper~[2] provides a detailed analysis of the main results on the formulation and solution of both classical and new modified [3] boundary value problems for the equation (\ref{eq:1}).

The main results on the formulation and study of boundary value problems for the equation (\ref{eq:1}) are obtained under the conditions imposed on the parameters of the equation: $0<|\alpha|, |\beta|, |\alpha+\beta| < 1$. Note that the Cauchy problem for the equation (\ref{eq:1}) at $\alpha=\beta$, $0< |\beta| < 1/2$ in the classical formulation with conditions
$$\lim_{y\to x+0} U(x,y) = \tau(x), \qquad  \lim_{y\to x+0} (y-x)^{2\beta} (U_y - U_x) = \nu(x), \quad y>x, $$
($\tau$, $\nu$ -- given functions) is invalid when $|\alpha| = |\beta| = \displaystyle\frac{1}{2}$ due to the fact that either the decision itself or its derivative in the normal direction (depending on the sign of $\beta$) on the line of singularity of the coefficients of $y=x$ goes to infinity.

In the present work, the authors propose to extend the formulation and solution of the modified Cauchy-type problems [4--6] for equation (\ref{eq:1}) in the following cases: a) $\alpha=\beta=\displaystyle\frac{1}{2}$; b) $\alpha=+\,\displaystyle\frac{1}{2}$, $\beta=-\,\displaystyle\frac{1}{2}$; c) $\alpha=\beta=-\,\displaystyle\frac{1}{2}$ (Problems $C_1$, $C_2$, $C_3$, respectively).

In the problem $C_1$ on the line $y=x$ conditions are set:
\begin{equation}
\lim_{y\to x+0} (y-x) (U_y - U_x) = \nu_1(x), \qquad 0\leq x < +\infty, \label{eq:3}
\end{equation}
\begin{equation}
\lim_{y\to x+0} \left[ U(x,y) - \nu_1(x) \left( \ln\sqrt{y-x} + \psi(\frac{1}{2}) - \psi(1)\right)\right] = \tau_1(x), \label{eq:4}
\end{equation}
in the problem $C_2$:
$$
\lim_{y\to x+0}  U(x,y) = \tau_2(x), \qquad 0\leq x < +\infty;
$$
$$
\lim_{y\to x+0} \left[ (U_y - U_x) - \frac{d}{dx} U(x,x) \left( \ln\sqrt{(y-x)} + 2\psi(\frac{1}{2}) - 2\psi(1)+1\right)\right] = \nu_2(x),$$
$$ \quad 0\leq x < +\infty,
$$
so in $C_3$:
$$
\lim_{y\to x+0}  U(x,y) = \tau_3(x), \qquad 0\leq x < +\infty;
$$
$$
\lim_{y\to x+0} \left[ (U_y - U_x)(y-x)^{-1} - \frac{d^2}{dx^2} U(x,x) \left(\ln\sqrt{(y-x)} + \psi(\frac{3}{2}) - \psi(3)\right)\right] = $$
$$ = \nu_3(x), \quad 0\leq x < +\infty,
$$

The problem $C_1$ is solved by the Riemann method, $C_2$, $C_3$ are obtained from the formula of the General solution of the Euler--Darboux equation.

On the basis of the obtained results the unambiguous solvability is proved
a modified $ \Delta_1^S$ problem in the domain representing the first quadrant, with boundary conditions on the coordinate axes and conjugation on the $y=x$ line.

\section*{$\Delta^S_1$ Problem} On the set ${\tt D}$ find the solution of the equation,
satisfying the conditions:
\begin{equation}
U(0,y) = \varphi_1(y), \qquad  0\leq y < +\infty, \label{eq:19}
\end{equation}
\begin{equation}
U(x,0) - \, \frac{1}{\Gamma^2(\frac{1}{2})} \int\limits^x_0 \nu_2(t) \, t^{- 1/2} (x-t)^{- 1/2} \ln\frac{t}{x} dt = \varphi_2(x), \qquad  0\leq x < +\infty. \label{eq:20}
\end{equation}

On the singularity line of the coefficients of the equation
the matching conditions
\begin{equation}
\nu_1(x) = \lim_{y\to x+0} (y-x) (U_y - U_x) = - \lim_{y\to x-0} (x-y) (U_x - U_y)= - \nu_2(x),  \label{eq:21}
\end{equation}
$$ \tau_1(x) = \lim_{y\to x+0} \left[ U(x,y) - \nu_1(x) \left( \ln\sqrt{y-x} + \psi(\frac{1}{2}) - \psi(1)\right)\right] =$$ \begin{equation}= \lim_{y\to x-0} \left[ U(x,y) - \nu_2(x) \left( \ln\sqrt{x-y} + \psi\left(\frac{1}{2}\right) - \psi(1)\right)\right] = \tau_2(x). \label{eq:22}
\end{equation}

Given functions $\varphi_k$, $k=1, 2$ are \ underline{conditions}: $\varphi_k(0)=0$, $\varphi_k(x) \in C^{(3)}[0, + \infty)$, $k=1, 2$.

The uniqueness of the solution of the $ \Delta^S_1$ problem follows from the uniqueness obtained by the Riemann method, the solution of the modified Cauchy problem taken as a basis, and the unique solvability of the system of integral equations obtained in the process of solving the problem. The existence of the solution is proved by verification.

\section*{Conclusion} Solutions of three modified Cauchy-type problems
for the Euler--Darboux equation in the case of parameters equal modulo one second were received. Based on the solution of the modified Cauchy problem, the unique solvability of generalizations of $\Delta_1$-problem with data on coordinate axes and conjugation conditions on the singularity line of the equation coefficients is proved. The solution of the problem $\Delta_1$ is obtained explicitly.

\vskip10pt

\medskip \centerline{\textsc{References}} \smallskip \small

\begin{enumerate}
\item {Smirnov M.\,M.,} {Degenerated Hyperbolic Equations, Vysshaya Skola, Minsk (1977).}
\item{Nakhushev A.\,M.,} {Certain boundary value problems for hyperbolic equations and equations of mixed type, Differ. Uravn., \textbf{5}:1 (1969), 44--59.}
\item{M. S. Salakhitdinov, B. I. Islomov, A nonlocal boundary-value problem with conormal derivative for a mixed-type equation with two inner degeneration lines and various orders of degeneracy, Izv. Vyssh. Uchebn. Zaved. Mat., 2011, no. \textbf{1}, 49--58; Russian Math. (Iz. VUZ), \textbf{55}:1 (2011), 42--49.}
\item {M. V. Dolgopolov, I. N. Rodionova, V. M. Dolgopolov, On one nonlocal problem for the Euler--Darboux equation, Vestn. Samar. Gos. Tekhn. Univ., Ser. Fiz.-Mat. Nauki [J. Samara State Tech. Univ., Ser. Phys. Math. Sci.],  \textbf{20}:2 (2016),  259--275.}
\item {I. N. Rodionova, V. M. Dolgopolov, M. V. Dolgopolov, Delta-problems for the generalized Euler--Darboux equation, Vestn. Samar. Gos. Tekhn. Univ., Ser. Fiz.-Mat. Nauki [J. Samara State Tech. Univ., Ser. Phys. Math. Sci.], \textbf{21}:3 (2017),  417--422.}
\item {Dolgopolov M.\,V., Dolgopolov V.\,M., Rodionova I.\,N., Two Problems Involving Equations of Hyperbolic Type of the Third Order in the Three-Dimensional Space, Advancement and Development in Mathematical Sciences, \textbf{3}:1--2 (2012), 25--38.}
\end{enumerate}
\end{document}